\documentclass[reqno]{amsart}
\usepackage{amsthm}
\usepackage{amsmath}
\usepackage{amssymb}
\usepackage{seqsplit}
\usepackage{url}

\newtheorem{theorem}{Theorem}[section]
\numberwithin{equation}{section}
\newcommand{\abs}[1]{\left\lvert#1\right\rvert}
\newcommand{\rd}[1]{\left\lfloor#1\right\rceil}

\begin{document}

\title{The size of oscillations in the Goldbach conjecture}

\author[M.~J. Mossinghoff]{Michael J. Mossinghoff}
\address{Center for Communications Research\\
Princeton, NJ, USA}
\email{m.mossinghoff@idaccr.org}

\author[T.~S. Trudgian]{Timothy S. Trudgian}\thanks{This work was supported by a Future Fellowship (FT160100094 to T.~S. Trudgian) from the Australian Research Council.}
\address{School of Science\\
UNSW Canberra at ADFA\\
ACT 2610, Australia}
\email{t.trudgian@adfa.edu.au}

\keywords{Goldbach conjecture, Hardy--Littlewood conjectures, oscillations, Riemann hypothesis, simultaneous approximation}
\subjclass[2010]{Primary: 11P32; Secondary: 11J20, 11M26, 11Y35}
\date{\today}

\begin{abstract}
\noindent
Let $R(n) = \sum_{a+b=n} \Lambda(a)\Lambda(b)$, where $\Lambda(\cdot)$ is the von Mangoldt function.
The function $R(n)$ is often studied in connection with Goldbach's conjecture.
On the Riemann hypothesis (RH) it is known that $\sum_{n\leq x} R(n) = x^2/2 - 4x^{3/2} G(x) + O(x^{1+\epsilon})$, where $G(x)=\Re \sum_{\gamma>0} \frac{x^{i\gamma}}{(\frac{1}{2} + i\gamma)(\frac{3}{2} + i\gamma)}$ and the sum is over the ordinates of the nontrivial zeros of the Riemann zeta function in the upper half-plane.
We prove (on RH) that each of the inequalities $G(x) < -0.02093$ and $G(x)> 0.02092$ hold infinitely often, and establish improved bounds under an assumption of linearly independence for zeros of the zeta function.
We also show that the bounds we obtain are very close to optimal.
\end{abstract}

\maketitle

\section{Introduction}\label{Intro}

Let $\Lambda(n)$ denote the von Mangoldt function, and define $R(n)$ by
\begin{equation}\label{cumin}
R(n) = \sum_{a+b=n} \Lambda(a) \Lambda(b),
\end{equation}
where the sum is over positive integers $a$ and $b$ that sum to $n$.
This function arises naturally in the study of Goldbach's problem: clearly $R(n)>0$ precisely when $n$ is the sum of two positive prime powers.
The use of the von Mangoldt function makes the problem more amenable to analysis, and Goldbach's conjecture would follow if it could be shown that $R(n)$ were sufficiently large at even integers $n>2$.
It is natural then to study the average value of $R(n)$.
It is known that 
\begin{equation}\label{coriander}
\sum_{n\leq x} R(n) = \frac{1}{2} x^{2} + O\left(x^{2}(\log x)^{-A}\right),
\end{equation}
unconditionally, for any positive constant $A$.
In a series of articles in 1991, Fujii obtained improvements on the error term in \eqref{coriander} that are conditional on the Riemann hypothesis.
In the first of this series, he established \cite{Fujii1991P1} that
\[
\sum_{n\leq x} R(n) = \frac{1}{2}x^{2} + O(x^{3/2}),
\]
and in the second paper \cite{Fujii1991P2} he refined the error term, proving that\footnote{We note that the sum over zeros, written in the form $2 \sum_{\rho}\frac{x^{\rho+1}}{\rho(\rho+1)}$ appears here even without assuming the Riemann hypothesis.}
\begin{equation}\label{ginger}
\sum_{n\leq x} R(n) = \frac{1}{2}x^{2} -4x^{3/2}\Re\sum_{\gamma>0} \frac{x^{i\gamma}}{(\frac{1}{2} + i\gamma)(\frac{3}{2} + i\gamma)} + O\left((x \log x)^{4/3}\right).
\end{equation}
Similar statements, with a slightly larger power on the $\log x$ term, were proved by Goldston \cite{Goldston} and by Granville \cite{G2}.
Reductions in the error term in \eqref{ginger} were made by Bhowmik and Schlage-Puchta \cite{BSP} and then by Languasco and Zaccagnini \cite{LZ}, who established $O(x \log^{3} x)$.
This is fairly close to optimal, since Bhowmik and Schlage-Puchta also proved that the error term here is $\Omega(x \log\log x)$. Analogous results for forms of \eqref{cumin}, where $n$ is written as the sum of $k$ prime powers, have been proved by Languasco and Zaccagnini \cite{LZ} and by Bhowmik, Ramar\'{e}, and Schlage-Puchta \cite{BRS}.

In this article, we study the oscillations in the sum on the right side of \eqref{ginger}.
To this end, define $G(x)$ by 
\begin{equation}\label{turmeric}
G(x) = \Re \sum_{\gamma>0} \frac{x^{i\gamma}}{(\frac{1}{2} + i\gamma)(\frac{3}{2} + i\gamma)},
\end{equation}
where the sum is over the ordinates of the zeros of the Riemann zeta function in the upper half-plane.
We assume the Riemann hypothesis, so each such zero has real part $1/2$.

In his third paper of 1991 on this topic \cite{Fujii1991P3}, Fujii proved that if the ordinates of the first $70$ zeros of the Riemann zeta function on the critical line are linearly independent over the rationals, then each of the inequalities
\begin{equation}\label{eqnFujiiBounds}
G(x) < -0.012,
\quad
G(x) > 0.012
\end{equation}
would hold for an unbounded sequence of positive real numbers $x$.
He noted that this conclusion could also be established without the linear independence hypothesis, if one instead employed a method of Odlyzko and te Riele to solve certain simultaneous approximation problems involving these $70$ real numbers.
In 1985 Odlyzko and te Riele \cite{OTR} famously employed this method to disprove the Mertens conjecture regarding the size of oscillations in the function $M(x) =\sum_{n\leq x} \mu(n)$, where $\mu(\cdot)$ represents the M\"obius function.
Recently, Hurst \cite{Hurst} used the same method, along with additional techniques, to obtain the presently best known result in this problem.

Odlyzko and te Riele established large oscillations in the positive direction by determining a real number $y$ and integers $m_1$, \ldots, $m_{70}$ with the property that
\[
\abs{\gamma_{k_j} y - \psi_{k_j} - 2m_j\pi} < \epsilon_1,
\]
for $1\leq j\leq70$, for a small positive number $\epsilon_1$.
Here $\psi_{k_j}$ represents the argument of the residue of $1/\zeta(s)$ at $s=1/2+i\gamma_{k_j}$, and $1\leq k_1<k_2<\cdots<k_{70}\leq400$ denotes a particular sequence of positive integers corresponding to the zeros which produced the most beneficial contributions in the method employed there.
Likewise, to establish large oscillations in the negative direction, they determined $z$, $n_1$, \ldots, $n_{70}$
so that
\[
\abs{\gamma_{k_j} z - \psi_{k_j} - (2n_j+1)\pi} < \epsilon_2,
\]
for $1\leq j\leq 70$, for a small positive number $\epsilon_2$.
In \cite{Fujii1991P3}, Fujii required analogous results for the same problems, but with each $\psi_{k_j}$ eliminated, $k_j=j$ for each $j$, and $\epsilon_1 = \epsilon_2 = 0.1$.
(The first case is then a homogeneous approximation problem, and one naturally also requires $y\neq0$ there.)
It is not clear however if the required computations were in fact performed in \cite{Fujii1991P3}: it is stated that the argument there implies the bounds \eqref{eqnFujiiBounds} ``in principle.''

In this article, we analyze the oscillations in $G(x)$, and prove two main results.
First, we use the method of Odlyzko and te Riele to establish a lower bound on the oscillations exhibited by this function in each direction, improving \eqref{eqnFujiiBounds}.
We also establish improved bounds under an assumption of linear independence for the zeros of the zeta function.
Second, we establish an upper bound on these oscillations, which shows that our results are close to optimal.
We prove the following theorem.

\begin{theorem}\label{thmMain}
With $G(x)$ as in \eqref{turmeric}, on the Riemann hypothesis each of the following inequalities holds for an unbounded sequence of positive real numbers $x$:
\begin{equation}\label{paprika}
G(x) < -0.020932, \quad G(x) > 0.020927.
\end{equation}
Moreover, for all $x>0$,
\begin{equation}\label{saffron}
\abs{G(x)} < 0.023059.
\end{equation}
In addition, if the ordinates of the first $10^6$ zeros of the Riemann zeta function in the upper half-plane are linearly independent over $\mathbb{Q}$, then each of the following inequalities holds for an unbounded sequence of positive real numbers $x$:
\begin{equation}\label{fenugreek}
G(x) < -0.022978, \quad G(x) > 0.022978.
\end{equation}
\end{theorem}

This paper is organized in the following way.
Section~\ref{sage} establishes the upper bound \eqref{saffron} of Theorem~\ref{thmMain}.
Section~\ref{secCondLB} obtains lower bounds for oscillations in $G(x)$, conditioned on the existence of solutions to particular simultaneous approximation problems involving a number of zeros of the Riemann zeta function, and establishes \eqref{fenugreek}. 
Last, Section~\ref{secComputations} describes the calculations required to establish the bounds \eqref{paprika} on the oscillations in this function without assuming any linear independence conditions to complete the proof of Theorem~\ref{thmMain}.

We remark that Hardy and Littlewood \cite{HL23} conjectured that $R(n) \sim nS(n)$ for even integers $n$, where
\begin{equation*}\label{cloves}
S(n) = \prod_{p|n} \left(1 + \frac{1}{p-1}\right) \prod_{p\nmid n} \left( 1 - \frac{1}{(p-1)^{2}}\right),
\end{equation*}
and that several authors encounter $G(x)$ when estimating the average value of $R(n) - n S(n)$.
For example, Fujii \cite{Fujii1991P2} in fact established \eqref{ginger} in the form
\[
\sum_{n\leq x} \bigl(R(n)-n S(n)\bigr) = -4 x^{3/2} G(x) + O\left((x\log x)^{4/3}\right).
\]
It is readily seen that the two forms are equivalent, since from Montgomery and Vaughan \cite[Lem.~1]{MV} we have that
\[
\sum_{n\leq x} n S(n) = \frac{1}{2} x^2 + O(x\log x).
\]
Additional estimates involving $S(n)$ and related functions and their application in problems in additive number theory can be found in \cite{MV}.

\section{An upper bound for $\abs{G(x)}$}\label{sage}

Taking the real part of the sum in \eqref{turmeric} produces
\begin{align}
\label{rosemary}
G(x) &= -\sum_{\gamma>0} \frac{\cos(\gamma \log x)}{\gamma^{2} +\frac{1}{4}}
+ \sum_{\gamma>0} \frac{3\cos(\gamma\log x) + 2\gamma \sin(\gamma \log x)}{(\gamma^{2} + \frac{1}{4})(\gamma^{2} + \frac{9}{4})}\\
\label{R1}
&= \sum_{\gamma>0} \frac{ (\frac{3}{4} - \gamma^{2}) \cos(\gamma \log x) + 2\gamma \sin(\gamma \log x)}{(\gamma^{2} + \frac{1}{4})(\gamma^{2} + \frac{9}{4})}.
\end{align}
With a little calculus one can show that the maximal value of the numerator in \eqref{R1} is $\sqrt{ \gamma^{4} + \frac{5}{2} \gamma^{2} + \frac{9}{16}}$, occurring when 
\[
\tan(\gamma\log x) = \frac{ 2\gamma}{\frac{3}{4} - \gamma^{2}},
\]
and that the minimal value is $-\sqrt{ \gamma^{4} + \frac{5}{2} \gamma^{2} + \frac{9}{16}}$, so
\begin{equation}\label{R4}
|G(x)| \leq \sum_{\gamma>0} h(\gamma), \quad h(\gamma) = \frac{\sqrt{ \gamma^{4} + \frac{5}{2} \gamma^{2} + \frac{9}{16}}}{(\gamma^{2}+\frac{1}{4})(\gamma^{2} + \frac{9}{4})}.
\end{equation}
A simple expansion shows that
\[
h(\gamma) = \frac{1}{\gamma^{2}+ \frac{1}{4}} - \frac{1}{\gamma^{4}} + \frac{2}{\gamma^{6}} - \frac{61}{16 \gamma^{8}} + O\left(\gamma^{-10}\right),
\]
and from Davenport \cite[ch.\ 12]{Davenport} we have that
\begin{equation}\label{porto}
\sum_{\gamma>0} \frac{1}{\gamma^{2} + \frac{1}{4}} = \sum_{\rho} \Re \left(\rho^{-1}\right) = 1 + \frac{\xi}{2} - \frac{\log 4\pi}{2} = 0.02309\ldots,
\end{equation}
where $\rho=1/2+i\gamma$ and $\xi=0.577\ldots$ represents Euler's constant.
We write
\begin{equation}\label{Q2}
h(\gamma) = \frac{1}{\gamma^{2} + \frac{1}{4}} - \frac{1}{\gamma^{4}} + U(\gamma).
\end{equation}
A simple calculation reveals that $U(\gamma)\gamma^{6} \leq 2$.
Therefore, to obtain an upper bound on $\abs{G(x)}$, we require an upper bound on $\sum_{\gamma>0} \gamma^{-6}$.
(We also need a lower bound on the sum over $\gamma^{-4}$ from \eqref{Q2}, but clearly any finite sum will work.)
For this, we employ the result of Lehman \cite[Lem.\ 3]{Lehman} stating that
\begin{equation}\label{Q2a}
\sum_{\gamma>T} \gamma^{-n} < \frac{\log T}{T^{n-1}}
\end{equation}
provided $T\geq 2\pi e=17.079\ldots$ and $n\geq 2$.
Using \eqref{R4}, \eqref{porto}, \eqref{Q2},  and \eqref{Q2a}, we therefore conclude that
\[
|G(x)| < 1 + \frac{\xi}{2} - \frac{\log 4\pi}{2} - \sum_{0< \gamma \leq T_1} \frac{1}{\gamma^{4}} + 2 \sum_{0<\gamma\leq T_2} \frac{1}{\gamma^{6}} + \frac{\log T_2}{T_2^{5}},
\]
where we may choose any values for $T_1>0$ and $T_2\geq2\pi e$.
Choosing the first $1000$ zeros for each sum, that is, taking $T_1 = T_2 = 1420.41$, we find that $|G(x)| < 0.023058681$, which establishes \eqref{saffron}.

\section{Conditional lower bounds}\label{secCondLB}

We may determine lower bounds on the oscillations of $G(x)$, conditioned on the existence of solutions to certain simultaneous approximation problems involving a number of nontrivial zeros of the Riemann zeta function.
We treat large oscillations in the positive direction here; large displacements in the negative direction follow analogously.

Given a positive integer $N$ and a positive real number $\epsilon$, suppose there exists a real number $y$ and integers $m_1$, \ldots, $m_N$ so that
\begin{equation}\label{oregano}
\abs{\gamma_k y - (2m_k+1)\pi} \leq \epsilon
\end{equation}
for $1\leq k\leq N$.
Then certainly
\begin{equation}\label{saltpepper}
\cos(\gamma_k y) < -1 + \frac{\epsilon^2}{2}
\end{equation}
for each $k$.
Let $T>2\pi e$ be a real number selected so that the number of nontrivial zeros of the Riemann zeta function with ordinate $\gamma<T$ is exactly $N$.
From \eqref{rosemary}, we have
\begin{equation}\label{cocoa}
\begin{split}
G(x) &= -\sum_{\gamma\leq T} \frac{\cos(\gamma \log x)}{\gamma^2+\frac{1}{4}}
+ \sum_{\gamma\leq T} \frac{3\cos(\gamma \log x)+2\gamma\sin(\gamma \log x)}{(\gamma^2+\frac{1}{4})(\gamma^2+\frac{9}{4})}\\
&\qquad + \sum_{\gamma>T} \frac{(\frac{3}{4}-\gamma^2)\cos(\gamma \log x) + 2\gamma\sin(\gamma \log x)}{(\gamma^2+\frac{1}{4})(\gamma^2+\frac{9}{4})}\\
&=: G_1(x,T) + G_2(x,T) + G_3(x,T).
\end{split}
\end{equation}
From \eqref{R4}, we have
\[
\abs{G_3(x,T)}
\leq\sum_{\gamma> T} \frac{\sqrt{\gamma^4+\frac{5}{2}\gamma^2+\frac{9}{16}}}{(\gamma^2+\frac{1}{4})(\gamma^2+\frac{9}{4})}
< \sum_{\gamma>T} \frac{1}{\gamma^2},
\]
and from \cite[Lem.\ 1]{Lehman} we obtain
\begin{equation}\label{mace}
\sum_{\gamma>T} \frac{1}{\gamma^2} = \frac{1}{2\pi}\int_T^\infty \frac{\log(t/2\pi)}{t^2}\,dt + \vartheta\left(\frac{4}{T^2}\log T + 2\int_T^\infty \frac{dt}{t^3}\right),
\end{equation}
where $\vartheta$ is a complex number satisfying $\abs{\vartheta}\leq1$.
This gives a better estimate than that in \eqref{Q2a}, which we shall need in what follows.
While the constants in the error in \eqref{mace} could be improved by the results in \cite{PT,TST}, the range of $T$ that we are considering here makes any potential gain negligible. 
Consequently,
\begin{equation}\label{nutmeg}
\abs{G_3(x,T)} <  B_3(T) :=  \frac{1}{2\pi T}\left(\log T + 1-\log 2\pi + \frac{2\pi}{T}(1+4\log T)\right)
\end{equation}
for all $x>0$.
For $G_1$, we use \eqref{saltpepper} to find
\[
G_1(e^y,T)  > \left(1-\frac{\epsilon^2}{2}\right)\sum_{\gamma\leq T} \frac{1}{\gamma^2+\frac{1}{4}} .
\]
For $G_2$, we observe that $3\cos t + 2\gamma\sin t$ is decreasing near $t=\pi$, so
\[
G_2(e^y,T) \geq -\sum_{\gamma\leq T} \frac{3\cos\epsilon+2\gamma\sin\epsilon}{(\gamma^2+\frac{1}{4})(\gamma^2+\frac{9}{4})}.
\]
Therefore,
\begin{equation}\label{eqnCondUpper}
G(e^y) > \left(1-\frac{\epsilon^2}{2}\right)\sum_{\gamma\leq T} \frac{1}{\gamma^2+\frac{1}{4}} - \sum_{\gamma\leq T} \frac{3\cos\epsilon+2\gamma\sin\epsilon}{(\gamma^2+\frac{1}{4})(\gamma^2+\frac{9}{4})} - B_3(T).
\end{equation}

Similarly, given $N$ and $\epsilon$, if we suppose there exists a real number $z$ and integers $m_1$, \ldots, $m_N$ so that
\begin{equation}\label{cardomom}
\abs{\gamma_k z - 2m_k\pi} \leq \epsilon
\end{equation}
for $1\leq k\leq N$, then we obtain the negation of the expression in \eqref{eqnCondUpper} as a lower bound on the oscillations of $G(x)$:
\[
G(e^z) < \left(\frac{\epsilon^2}{2}-1\right)\sum_{\gamma\leq T} \frac{1}{\gamma^2+\frac{1}{4}} + \sum_{\gamma\leq T} \frac{3\cos\epsilon+2\gamma\sin\epsilon}{(\gamma^2+\frac{1}{4})(\gamma^2+\frac{9}{4})} + B_3(T).
\]

In Table~\ref{tableCondEst} we list a few values for the bound \eqref{eqnCondUpper} for a number of choices of $N$.
In each case we assume $\epsilon=0.01$, and take $T = T^*(N)$, where
\begin{equation}\label{allspice}
T^*(N) =  \gamma_{N+1} - \frac{\gamma_{N+1}-\gamma_N}{100}.
\end{equation}

\begin{table}[tbh]
\caption{Conditional lower bounds for large positive values of $G(x)$ from \eqref{eqnCondUpper}, assuming the simultaneous approximation problem \eqref{oregano} has a solution with $\epsilon=0.01$.}\label{tableCondEst}
\begin{tabular}{|cc|cc|}\hline
$N$ & Bound\ & $N$ & Bound\\\hline
  70 & $0.014756$ & 500 & $0.020630$\\
100 & $0.016352$ & 600 & $0.020902$\\
150 & $0.017837$ & 700 & $0.021109$\\
200 & $0.018692$ & 800 & $0.021272$\\
250 & $0.019269$ & 900 & $0.021404$\\
300 & $0.019684$ & 1000 & $0.021515$\\
350 & $0.020001$ & 2000 & $0.022079$\\
400 & $0.020254$ & $10^4$ & $0.022699$\\
450 & $0.020459$ & $10^5$ & $0.022925$\\\hline
\end{tabular}
\end{table}

If the ordinates of the first $N$ nontrivial zeros of the zeta function are linearly independent, then by Kronecker's theorem the corresponding bound in Table~\ref{tableCondEst} would necessarily follow, as would any value computed with an arbitrary choice of $\epsilon>0$.
Selecting $\epsilon=10^{-6}$ with $N=10^6$ produces the value $0.02297864\ldots$\,, which verifies \eqref{fenugreek} in Theorem~\ref{thmMain}.

To obtain bounds without linear independence, in the next section we turn to the method of Odlyzko and te Riele for constructing solutions to some of these simultaneous approximation problems.

\section{Computations}\label{secComputations}

We complete the proof of Theorem~\ref{thmMain} by solving the simultaneous approximation problems \eqref{oregano} and \eqref{cardomom} for particular $N$ and $\epsilon$.
For this we employ the method of Odlyzko and te Riele \cite{OTR}, which we briefly describe here.
Let $\rd{x}$ denote the integer nearest the real number $x$, and let $\mathbf{e}_k$ denote the $k$th elementary unit column vector in the appropriate real vector space.
The construction requires values for four integer parameters: $N$, $b$, $c$, and $d$.
Here, $b$ represents the number of bits of precision used in the computation; $c$ and $d$ are small positive integers whose meanings will be described shortly.

\subsection{Large positive values for $G(x)$}\label{subsecHigh}
Consider first the inhomogeneous problem \eqref{oregano}, where we require a real number $y$ with the property that $\gamma_k y$ is near $\pi$, modulo integer multiples of $2\pi$, for $1\leq k\leq N$.
We construct the $(N+2)\times(N+2)$ integer matrix $M$ whose column vectors are
\begin{gather*}
\rd{2^{b+1}\pi}\mathbf{e}_k,\; 1\leq k\leq N,\\
\mathbf{e}_{N+1}-\sum_{k=1}^N \rd{2^{b-c} \gamma_k}\mathbf{e}_k,\\
2^b N^d \mathbf{e}_{N+2}+\rd{2^b\pi}\sum_{k=1}^N \mathbf{e}_k.
\end{gather*}
That is, $M$ consists of an $(N+2)\times N$ diagonal matrix with entries $\rd{2^{b+1}\pi}$ on the diagonal, augmented with one column carrying rounded multiples of the $\gamma_k$, and another largely filled with a rounded multiple of the inhomogeneous part, $\pi$.
The penultimate vector carries the lone nonzero value in vector position $N+1$, set to $1$ so that we can recover a coefficient later in the computation.
The last vector has the only nonzero value in the last position, chosen to be much larger than the other entries of the matrix.

We apply the LLL algorithm \cite{LLL} to $M$ to compute a reduced basis for the lattice spanned by its column vectors.
This reduced basis consists of vectors that are relatively short, in fact within a factor (whose value is bounded by an expression that is exponential in the dimension) of the shortest independent vectors in the lattice.
Since the last coordinate of every vector in the lattice is an integer multiple of the large integer $2^b N^d$, it is likely that there is only one vector in the reduced lattice with a nonzero value in this position, which is very likely to be $\pm2^b N^d$.
If this value is negative we can negate the vector, so suppose it is $(r_1,\ldots,r_N,s,2^b N^d)^T$.
We then have that there exist integers $m_1$, \ldots, $m_N$ such that
\[
r_k = m_k \rd{2^{b+1}\pi} + \rd{2^b\pi} - s\rd{2^{b-c}\gamma_k}
\]
for $1\leq k\leq N$, and that the $r_k$ are relatively small.
If $s<0$ then we can negate this vector so that our inhomogeneous part is $-\pi$, which serves us just as well, so we assume $s\geq0$ here.
We might then expect
\[
\gamma_k s 2^{-c} \approx 2\pi  m_k+ \pi
\]
so we take $y=s/2^c$, and use this in \eqref{cocoa} and \eqref{nutmeg} to compute the resulting lower bound on positive values reached by $G(x)$:
\begin{equation}\label{basil}
G_1(e^y,T^*(N)) + G_2(e^y,T^*(N)) - B_3(T^*(N)),
\end{equation}
with $T^*(N)$ as in \eqref{allspice}.
For each $k$ we also compute $m_k = \rd{(\gamma_k y - \pi)/2\pi}$, and then
\begin{equation}\label{eqnEps1}
\epsilon_1 = \max_{1\leq k\leq N}\{\abs{\gamma_k y - (2m_k+1)\pi}\}.
\end{equation}
A large value of $\epsilon_1$ (and consequently a small value in \eqref{basil}) likely indicates that insufficient precision was employed.
In that case we repeat this process with a larger value of $b$.

Odlyzko and te Riele used $c=10$ and $d=4$.
Both values worked sufficiently well in our application, too, so we did not alter these in our principal runs.
Those authors also reported selecting $b$ between $6.6N$ and $13.3N$ (that is, using between $2N$ and $4N$ decimal digits of precision).
The larger end of this range sufficed in our application only for $N$ up to about $250$, where we produced $\epsilon_1=0.035$.
For larger dimensions we needed to select $b$ as large as $25N$.

\subsection{Large negative values for $G(x)$}\label{subsecLow}
For this case, we need to solve the homogeneous simultaneous approximation problem \eqref{cardomom}, as we need to find a value $z$ so that $\gamma_k z$ is very near an integer multiple of $2\pi$, for each $k$.
No additional computations are required here, as our reduced basis from the prior computation already contains many vectors of the form $(r_1,\ldots,r_N,t,0)$, in fact, there are likely to be $N+1$ of these.
Each one represents a viable solution to the homogeneous problem, since here for each $k$ we have
\[
r_k = m_k \rd{2^{b+1}\pi} - t\rd{2^{b-c}\gamma_k}
\]
for some integer $m_k$, and again the $r_k$ value is relatively small, so we might expect
\[
\gamma_k t 2^{-c} \approx 2\pi m_k
\]
for each $k$.
We may assume $t\geq0$.
For each such vector, we set $z=t/2^c$, and compute the resulting bound on negative values achieved by $G(x)$:
\[
G_1(e^z,T^*(N)) + G_2(e^z,T^*(N)) + B_3(T^*(N)).
\]
Among all such vectors we may select the $z$ that produces the best value.
With this set, then for each $k$ we may compute $m_k = \rd{\gamma_k z/2\pi}$, and then
\begin{equation}\label{eqnEps2}
\epsilon_2 = \max_{1\leq k\leq N}\{\abs{\gamma_k z - 2m_k\pi}\}.
\end{equation}
Again, if $\epsilon_2$ is too large, then we can repeat the process with a larger value of $b$.
In practice, if $\epsilon_1$ was sufficiently small then $\epsilon_2$ was as well.

\subsection{Results}\label{subsecResults}

All computations were performed in SageMath \cite{Sage}, using resources at NCI Australia and at the Center for Communications Research.
High-precision values for zeros of the Riemann zeta function were computed using the \texttt{mpmath} Python library \cite{mpmath}, available within SageMath.

Table~\ref{tableResults} records the bounds we obtained on $G(x)$ in this way, using different values for $N$ and $b$.
The last line in this table records the parameters and results of the computation that establishes \eqref{paprika} in Theorem~\ref{thmMain}.
This calculation required almost two weeks of core time on an Intel Xeon Platinum 8175M processor running at $2.5$ GHz.
Figure~\ref{figYZ} exhibits the values of $2^{10}y$ and $2^{10}z$ obtained for this case, using base $36$ for economy of space.

\begin{table}[tbh]
\caption{Guaranteed oscillations in $G(x)$, along with the errors $\epsilon_1$ and $\epsilon_2$ from \eqref{eqnEps1} and \eqref{eqnEps2}, obtained by solving the simultaneous approximation problems using the first $N$ zeros of the Riemann zeta function, and using $b$ bits of precision. The displayed values for the bounds are truncated at the last displayed digit; those for $\epsilon_1$ and $\epsilon_2$ are rounded up at the last displayed digit.}\label{tableResults}
\begin{tabular}{|cccccc|}\hline
$N$ & $b$ & Lower & Upper & $\epsilon_2$ & $\epsilon_1$\\\hline
70 & 930 & $-0.0147727$ & $0.0147720$ & $0.00092$ & $0.00089$\\
100 &  1330 & $-0.0163668$ & $0.0163683$ & $0.00168$ & $0.00125$\\
150 &  2000 & $-0.0178557$ & $0.0178520$ & $0.00444$ & $0.00394$\\
200 &  2660 & $-0.0186992$ & $0.0187115$ & $0.01340$ & $0.01160$\\
250 &  3325 & $-0.0192583$ & $0.0192902$ & $0.03166$ & $0.03525$\\
300 &  4500 & $-0.0197172$ & $0.0196887$ & $0.02851$ & $0.03691$\\
350 &  6000 & $-0.0200320$ & $0.0200230$ & $0.02075$ & $0.01719$\\
400 &  7600 & $-0.0202570$ & $0.0202690$ & $0.01772$ & $0.00990$\\
450 &  8600 & $-0.0204387$ & $0.0204629$ & $0.04154$ & $0.03556$\\
500 & 11000 & $-0.0206646$ & $0.0206304$ & $0.02424$ & $0.02210$\\
600 & 15000 & $-0.0209324$ & $0.0209272$ & $0.02479$ & $0.02106$\\\hline
\end{tabular}
\end{table}

\begin{figure}[tbh]
\caption{$2^{10}y$ and $2^{10}z$ (in base $36$) for the last line in Table~\ref{tableResults}.}\label{figYZ}
\begin{tabular}{p{4.85in}}
\tiny\seqsplit{%
86frdmdbt0zjkhp0qj2lepxq6ei4fv1nz2bhuosu62e0z2jt9sgludvxmlbc8o8h4jccdgp6e9le5dqzv8jh9je13bd76p4cwady53dzyxnvwvv5llvdcznlk0owiwmzx1qigyn3ahiuqm4r4jiktbujq877u2osx6yvzzhqjrkcayo7095yenqhwckngl1wsj1pj93c4p2e0m2rzsia2ikdfmurdbzjcrm5gk89mas4lfo7qlbaoaqkh7sceaesdv6ghxtzaetf2x0zrom7bzm01io89k6lvwas64wz6xsebhzxa5l5r7bjvjmt3ly7px85qi94n70tgbtf8y97mvc7ie907pi2eto9759k60s3clald9nzk2i6mc1nbif245b64p2i5z01lr96wghldyl9q2gz3b4om8hq044u30afsbhxg3850fmjc06rdikvsgsojhzldkb9v7s8hw6ylyp12hn2bdxslvhhnv099qyma8bxiatdmwwjeefbryxyyf5g1lfih0ci74nxepotixswukq4t2y7u8i2bf1y4v1b0w0tdhvr6wmskat5jdhqjruohvn4zel9tf0fu2ensnc3b7jszs84h9k1ylk9pw31hfyrkk0unio39pbd4ntl50vlfgurx6h7ts0he7ofz1jtuz5cw11fuhotse36kcxro9tgb7ejaxc3azxt4kv63kbbbc4a621jzdnjm3oconqpv01m2uj4sbrglk9kdadquosarp6kwwea5rn9lqgd68kvv9rdncmientrp09w9z8i10pvjabeow6u3mxbecxjroz4jgs65z6m9c20feaf9rimj5mw5eokvcblqhojt2hrybqgwyvx842rwnez269i1dmdattnf8a9wbu2qrhz9l86qngccqefgxvduvhnojzdiqtwkwnkrj87xf5pnsc3or2zvo107flxpc4q0h6yhanr4vusc7y531oinauwer99oxkt8wkbg4tx298uycoiirua5wcfvp97zkyc8p8dnwjclpxg48otil6lk6hrjkr1vo17jxab4ftuqqtx05k38kch4v0fkq8o4ex2bqp6akusvnfbq8fyapiraftnngxkigj63pms3fdlexml7hjz34pijedagsaelq4oychhsqqgx2j97hxch0tpgw0k09xu8xvaig565egy6jzk56erxbyupz06rvtjhbhm1ekwa5emjvqaexp429wnzrlvxid5us4s9fk0vffxe4ec1hzblh2d5rpnrv78uw5zzcip3ksib2jydrm71amletabf2phqn5lnfyj9cr9ajpoq024pfwxi5295g0y8ogms5r6v9djyn8739fgibnhhknkel6rm3els7z1bp2owrej9x8dslg3yy4tdvcnykkm97cbz7kbxjge4ceiig9x5rwio2z1xqd2bxbegj8xds1l1cz7le0msrz85gjlvk7usohyi3uxvlbpbhlnwwt0wxbk7w8xklpgeb3d8hkyullzcs3cjjv9owjqbvrr3lo79yrkcuo4bk1uuawfbgp7ewqg73otkupzzfpbmecnbdok5sawzg4zzhhyra6q7susf8x9odq0vs6u8iuwbgdbwu35blafwkp50y9y1m7kwmpzouln0bcbe8dqyembgluqkon3fhaznk91n4l12yeoknxmiebl9irei07rwapcbj9n9jbjch1k9gbm9dfapf05wm9mnyiq6ricmvg1nq9dnndicyeb529shs0kf09548buag6ktgc8d8ebnpzuf9ryl9f2jin6wsom2n2ntj0me0zpm50doswtqjkpj9a9rxvk2lh85ofo5tu64p6o61t23py2qnx5zeix5ptt8srxl0wybgue9lvqud3u9s5wsr78xtvbt1iloiqvq7e0v5v89o3vftq58cj6ml0fqgr8ppxl6tgx2i37aj16127xixcx9s4cmu8qq1o3xqzfbdvh37penpt1j30ij13b5574kxdr2rb2xcmikwe1sd6t4po75nkdpj7jlcwql0qg90l2vaynw8cb0vw9cnd3otmrwjif7oi5budpsjtu2iiup4tmskr4jfi8y9eoq369wjy163dcqfgs39gu5ogewgvq7teceecsm581a2w2be0an9pq1krpi3ahkyhp0wkuf6uh57vlbbc6bdzkclo34209f0598jzh4ukdygdzipof8fp47zvsnizfhhsu81lkl1ec2ctg6xg63q3fu25eidruznxyz25cad2a56f150vh9hg39kiov1wxq7j60qgxyd1dz8jprokp2cu2eit90rzgp0k5as2icx5rkb4oscxkb0pjok6md03bekfhwp8wnoxmvp616gfmoak13x4bjty5gcxrsd7tp3576t3qmg5la1lpjfa9ykih06hyzgfkwnc5z0tsl6dfxniy6qih1r82806s9jrc3z0q6bxytje3y0lsj9b8q7f0o8s1zfxhxqrzlm4pcavefgxxg7hbf70zuh8ewqot2eowyrapoccv4w3cpk661ymn0l7zmba7i1bclle5c0cb03l1h5ngyrll0bp60etf42d47ys20a4hi2ccbhbn0fikfqa26ri685fldc62pg2uy1o18kdsxwctik460k7f0sg9fovhhsqi88ly0iqsjjr98qmpyzfvw16m3bj9v6ccact9jn07azk51s6lp4csyd0zwn1ookto81rqre7rqigaxenx65vjlaxkqq862f2m55n74q8k1asfauf7t5etqanzefis080rabw3l5rutobx3vriri9hvn3yaugi4ccqnhkz3kjoc6muu1gyuziopnq6ewy5w6rz6gww3fhw2a2sguzv0i629wnirddbvkenjwl99jr2j%
}\\
\tiny\seqsplit{%
8e4o4jhjivnaqrka615olndxdsaxzd0r04efnwwtqv6zfd25awtt7g37mz1orzafgutr4gps99mrngrjh7suujgcukden9tzrvucbvllqcnn9160ssjyylx1jzr3wes0u9w9nqztz0vm5hq8xgk5mdzj8pzdqe1h48l6fbwnbffonq36dex61o225rh7ro7er3o2gqit4bta8umnmwdld1k8qfr6njtovjg5o7ar6tci5o7n66q6zpxmiszxc5jnyv5olf4capgireuycqu1wqsqqzqj54bvvczu89m4n7lu26n4whla55eo5jyddiglaurlchlrl0m6s99xdvw9aa3lczk85oae846kxv665a7shd32oyuvdwai20sn1y44sgwjc01id5ogxnf4ri3xtbgxd3gijpb5ixrqwz10brqock3xjvcmjzl66j00llqig0evjq828n4oir93dhpudyiq7c6zd867u12vd3bs4u53nsou57hj9jeg8hty0js03nziwrl0nce21a1g853k3q9fy3arvuduryp71tx7msq2sjd65fi22iqvh3va3ir3uef4tx0ce836ldy9k0ifvrrz4zf2t8nxel8jfdujlqqfhadnk8d5uvov8os7c4w2tijg4run6eyn38st56uyfth51ix8kdukuqdv7puhfhc1tnk1d0okvjdnft9ajbyzpje1fy1692kpksowscoe7oog8fl4xc85gar8p692jsk71ihpngj2bsamy29oncn5037jupc8j7xquqyrp4c3yicw0f3u8w7wi8a25ny8vhvbqx8z1x0tgjct2stfxb4gopot9goxgrqt2nertqtfwv982nqwhfb63iqr4b3rkgm2jg5bqbl7d25titvg03wm2eyt0am7vtcdfu4c2qa4pdwrgwjn4y03hc2kznoyqgfajalf5ykgdbskb1pnzk2kd3dmfm80ssji70mluu9q69rk4bmy0ao7s09knlemyuo900hazhpl8hd6744kgfxsf05xdtgur9qiyd0hd6gyu3kmgx27lyy60y3r8trzts6cq9br7dsosard7g980crligdvt2akyxbqsi6288p9ai6yugy5ab5pmgl7c5jzc1t1slr0gtbb0mc0utsrz4wif6jwr20lpr0szodh3lttgndbp914sv4fo6gy4t8wuii5jyfucp2avwyfaicsoforz5wgtxwjssq8wypnvq7c4ig6d4ny9yt3z7h1p8zt62gjnp6tpci077zg3ap14eya7adnrc2vqdtysoo18o6w8dspliphw4p96878m4kru4nesop3ne650ut7zabd3tc5xvpzcss132ojmfpmoc8yn02v8kah2cvyhftio43qarflz2ohxuo64c92iyvt8hyaehjnsn9df8h90c5djvyqfv21zsr2mexhb37r9vx9qacbkxemt05jitf6e7fs7ajo6k59j8gdanq5huuwcle5xp68p1crezlcyt2a8g8wexrsfznu82q6eau9p5jhi481y708nkc5rsor6ydc9hqs48r8uctuk6nmtl2mgop9ts3fllxjbqnyx8at4uqsp50s5a8whdgyxhkroeyougw4jppblrml9prthgfzp0w3d1nc61ic6azzzdue4yqfmwnpmefrwbum1pfp33rhtujtq8uwaxnf50i5gj430xz3jw6ke3xqx727g3ocyfy9em1yrqgi5ntyauypaa1j3cte0l3c03o5ybcj13x0tmoyggnvihchpjmxuqbgcblwzazzxwrq7hejtva85jk425ro5m7xt6tesbgpfa1mlhc626dv0xsps6qq6c1mornpkswoe4wsh8ivda8gmt156s841wvwc2u8rxm0iezmnict00cvaljfy3qgxrvpb4mztg0yhb8x67bfh35301vmme4gm1dbf31xdyc5irqr2p69y6ukpjhb10yliieb74d0992h3jm8r63249pstte29ml9bdlaybqy48uomr8filjkgff56vemnaqaziv9kk7w5svqv3h0d6zfyk82od9qrflvow0yduuhhl1v5a47aenw99xtkf42tklyfcxokgqs3fsrl3v40fudl71w2f6jflub5665hajb3b5w0vu71ocsrymf1ay8bu86cpb037mlrxr0uqwn57vltl0358pmpmpro38anr3az2dr6ndaeqkgf0gh6w33q10n50yypzmf3awu3qwbs6z8odhaioax51m1p3hsluc92nlrt6k4ygurctf3lkunyl2narewdrt6k3ixb8wd4lysp3l2gm2nxjslp2e2v8qyxstsxmr1t7yr5kvf66r7jnf9pxdp1ie87k4abuv5an8v3qz3gsyt8q8pvq1isauizhlgvf08xmfh2r8e1eqtnm4wbg6u3oe43qhohks3ssl7q7yliot6d1fdrymiquj5bw0tu7uau645ukdupu7wbojd3xbcml65vy6ji3x8dx5i951jadg5szy3y5sq1vhbqpcsasbj9ryketz0ef3tg0aa3k0zgbx68e13capii8tyxmbeuc2yzloek3qd5htv17jrynp3g0syx82uipllqeu1ce2qd0y5kg8441jwh4qfli2m4hed9i4k2klf7jmht18ub3too29g9rz0n24n541n4xurx999i695l3snqn8rhwwamh8p2fh1sumh8nxtq8nmg2oo20h7br88tqtl63idfev4enspeg305lsddtk8zb51msj7ymnw8djvntjod9grja06rdkypuib9ziwz5qmn0dtpmmyqhwxtwltpczn1glfy3eexo0hxae26pte5l6iblb9ebfvhcmzj16rbing5z47lvguneaydygtmgixo2%
}
\end{tabular}
\end{figure}

\section*{Acknowledgments}

We thank NCI Australia, UNSW Canberra, and the Center for Communications Research for computational resources.
This research was undertaken with the assistance of resources and services from the National Computational Infrastructure (NCI), which is supported by the Australian Government.

\bibliographystyle{amsplain}

\end{document}